\documentclass{amsart}
\usepackage{amsaddr}
\usepackage{hyperref}
\usepackage{etoolbox}
\usepackage{lipsum}

\usepackage{graphicx}
\usepackage{tikz}
\usepackage{tikz-cd}
\usepackage{verbatim}
\usepackage{setspace}
\usepackage{mathrsfs}
\usepackage{hyperref}
\usepackage{academicons}
\usepackage{xcolor}
\newtheorem{theorem}{Theorem}[section]
\newtheorem{lemma}[theorem]{Lemma}

\theoremstyle{definition}

\theoremstyle{remark}

\theoremstyle{proposition}
\newtheorem{proposition}[theorem]{Proposition}

\theoremstyle{corollary}
\newtheorem{corollary}[theorem]{Corollary}

\numberwithin{equation}{section}

\DeclareMathOperator{\Hom}{Hom}

\DeclareMathOperator{\E}{E}

\newcommand*\xbar[1]{%
\hbox{%
\vbox{%
\hrule height 0.5pt % The actual bar
\kern0.5ex% % Distance between bar and symbol
\hbox{%
\kern-0.1em% % Shortening on the left side
\ensuremath{#1}%
\kern-0.1em% % Shortening on the right side
}%
}%
}%
}
\newcommand\restr[2]{{% we make the whole thing an ordinary symbol
\left.\kern-\nulldelimiterspace % automatically resize the bar with \right
#1 % the function
\right|_{#2} % this is the delimiter
}}

\newcommand{\orcid}[1]{\href{https://orcid.org/#1}{\textcolor[HTML]{A6CE39}{\aiOrcid}}}
\begin{document}
%\begin{titlepage}
\title{Isomorphisms Between Injective Modules}

%    Information for first author
\author{Mohanad Farhan Hamid} %\orcid{0000-0002-5384-0148}}
%\title[Running title]{Title}
%\UniCountry{University of Zagreb, Croatia}
%\author[F. Author, S. Author and L. Author]{First author, Second author and Last author}
%\address[First author]{Department of Mathematics\\University of Zagreb\\ 10 000 Zagreb\\ Croatia}
%\email[First author]{glasnikFirst@math.hr}
\address{Department of Applied Siences,
\\University of Technology-Iraq, Baghdad, Iraq
%\\Postcode (ZIP) :10066
%\\Telephone: {+964 783 8 972 842}
\\E-mails: {mohanad.f.hamid@uotechnology.edu.iq
\newline mohanadfhamid@yahoo.com}}
%\item \href{https://orcid.org/0000-0002-5384-0148}{\textcolor{orcidlogocol}{\aiOrcid} \hspace{2mm} orcid.org/0000-0002-5384-0148}

 \curraddr{}

%    \thanks will become a 1st page footnote.
\thanks{}

%    Information for second author

%    Current address

%    General info
 %\keywords{(Self) pure submodule,
%absolutely (self) pure module, quasi injective module}
%\end{titlepage}

%\newpage

\maketitle

%%%%%%%%%%%%%%%%%%%%%%%%%%%%%%%%%%%%%%%%%%%%%%%%%%%%%%%%%%%%%%%%%%%%%%%%
%%%%%%%%%%%%%%%%%%%%%%%%%%%%%%%%%%%%%%%%%%%%%%%%%%%%%%%%%%%%%%%%%%%%%%%%
%\newpage
%\title{\bf{Absolutely Self Pure Modules}}
%\begin{abstract}
\section*{\bf{Abstract}}
\noindent Suppose that $(\mathcal{F},\mathcal{M})$ is an injective structure of $R$-Mod such that the class $\mathcal{F}$ is closed for direct limits, then two modules in $\mathcal{M}$ are isomorphic if there are maps in $\mathcal{F}$ from each one of the modules into the other. Examples of module classes in such injective structures include (pure, coneat, and RD-) injective modules, as well as $\tau$-injective modules for a hereditary torsion theory $\tau$. Thus providing a generalization of a classical result of Bumby's and two recent ones by Mac\'{i}as-D\'{i}az. 
%\noindent Suppose that $(\mathcal{F},\mathcal{M})$ is an injective structure of $R$-Mod such that the class $\mathcal{F}$ is closed for direct limits, then the class $\mathcal{M}$ must be enveloping. Moreover, any two modules  have isomorphic $\mathcal{M}$-envelopes whenever there are maps in $\mathcal{F}$ from each one of the modules into the other. Examples of module classes in such injective structures include (pure, coneat, and RD-) injective modules, as well as $\tau$-injective modules for a hereditary torsion theory $\tau$. Thus providing a generalization of a classical result of Bumby's and two recent ones by Mac\'{i}as-D\'{i}az. 

\bigskip
\noindent \textbf{\keywordsname.} Injective structure, injective module, coneat injective module, pure injective module, RD-injective module.

\section{\bf{Introduction}}
In 1965, R. T. Bumby \cite{bumby} proved that two injective modules isomorphic to submodules of each other must be isomorphic. J. E. Mac\'{i}as-D\'{i}az \cite{diaz1} gave an analogous result that, over a commutative domain, two pure (resp., RD-) injective modules isomorphic to pure (resp. RD-) submodules of each other must be isomorphic. Recently, he proved that, over a commutative ring, two coneat injective modules isomorphic to coneat submodules of each other must be isomorphic \cite{diaz2}. There is a unifying theme for those results and the idea of injective structures provides such a unification. Each of the above classes is an example of a module class in an injective structure, in the sense of Maranda \cite{maranda}. The isomorphisms, in the above results, between a (generalized) injective module and some type of submodule of another (generalized) injective module are members in the class of homomorphisms of the corresponding injective structure. So our main result states that for an injective structure ($\mathcal{F}$, $\mathcal{M}$), any two modules in the class $\mathcal{M}$ are isomorphic if and only if there are homomorphisms in the class $\mathcal{F}$ from each of the modules into the other.

%$\mathcal{PBI} =$ the class of pure Baer injective $R$-modules.

%$\mathcal{RCI} =$ the class of $R$-coneat injective $R$-modules.

 %the class of $\tau$-injective $R$-modules for some hereditary torsion theory $\tau$.
\section{\bf{Preliminaries}} 
Let $R$ be a ring possessing an identity. Modules, unless stated otherwise, are all left unitary $R$-modules. Let $\mathcal{F}$ be a class of homomorphisms between modules and $\mathcal{M}$ a class of modules. An $\mathcal{M}$-\emph{preenvelope} of a module $A$ \cite[p. 129]{eandj} is a map $A \rightarrow M$ into a module $M \in \mathcal{M}$, with respect to which every map $A \rightarrow M'$, with $M' \in \mathcal{M}$, is injective. That is, there is a map $f: M \rightarrow M'$ which, when composed with $A \rightarrow M$, gives $A \rightarrow M'$. An $\mathcal{M}$-\emph{envelope} of $A$ is an $\mathcal{M}$-preenvelope $A \rightarrow M$ such that the only endomorphisms of $M$ which, when composed with $A \rightarrow M$, give $A \rightarrow M$ again, are automorphisms of $M$. Thus an $\mathcal{M}$-envelope of a module $A$, if existing, is unique up to isomorphism and is a direct summand of any $\mathcal{M}$-preenvelope of $A$ \cite[p. 129]{eandj}. The class $\mathcal{M}$ is (\emph{pre})\emph{enveloping} if every module has an $\mathcal{M}$-(pre)envelope.

Following \cite{maranda}, $\Phi(\mathcal{F})$ denotes the class of all modules injective with respect to maps in $\mathcal{F}$, and $\Psi(\mathcal{M})$ denotes the class of all homomorphisms of $R$-modules with respect to which all modules in $\mathcal{M}$ are injective. 
A pair $(\mathcal{F},\mathcal{M})$ is called \emph{injective structure} of $R$-Mod \cite{maranda} if the following three conditions hold:
\begin{enumerate}
    \item $\Psi(\mathcal{M}) = \mathcal{F}$.
    \item $\Phi(\mathcal{F}) = \mathcal{M}$.
    \item $\mathcal{M}$ is preenveloping.
\end{enumerate}
Modules in $\mathcal{M}$ will be called $\mathcal{F}$-injective modules, and maps in $\mathcal{F}$ will be called $\mathcal{F}$-maps (or $\mathcal{F}$-homomorphisms). If the $\mathcal{F}$-map is also a monomorphism (resp., epimorphism) we say that it is an $\mathcal{F}$-monomorphism (resp., $\mathcal{F}$-epimorphism). In \cite[p. 322]{wis}, $\mathcal{F}$-monomorphisms are called $\mathcal{M}$-\emph{copure} monomorphisms. 

Homomorphisms in $\mathcal{F}$ and modules in $\mathcal{M}$ have the following properties. They are all straightforward or can be found in \cite{maranda}.
\begin{enumerate}
    \item Any direct summand of an $\mathcal{F}$-injective module is $\mathcal{F}$-injective.
    \item If $f$ and $g$ are $\mathcal{F}$-maps then so is $gf$, if defined.
    \item If $gf$ is an $\mathcal{F}$-map then so is $f$.
    \item Split maps are in $\mathcal{F}$.
\end{enumerate}

An example of an injective structure is when $\mathcal{M}$ is the class of all $R$-modules and therefore, $\mathcal{F}$ consists of split monomorphisms.
Another example is $(\Psi (\mathcal{I}nj), \mathcal{I}nj)$, where $\mathcal{I}nj$ denotes the class of injective $R$-modules. In this case $\Psi (\mathcal{I}nj)$ consists of all monomorphisms of $R$-modules. 

%By an $\mathcal{F}$-injective module we mean a module in $\mathcal{M}$, and by an $\mathcal{F}$-map (or $\mathcal{F}$-homomorphism) we mean a map in $\mathcal{F}$. If the $\mathcal{F}$-map is also a monomorphism (resp., epimorphism) we say that it is an $\mathcal{F}$-monomorphism (resp., $\mathcal{F}$-epimorphism).

An monomorphism $A\rightarrow B$ of $R$-modules is a \emph{pure} (resp., \emph{coneat}) monomorphism if $A \otimes U \rightarrow B \otimes U$ is again a monomorpism for every (simple) right $R$-module $U$ \cite[p. 286]{wis} (\cite{scrivei}). It is called a \emph{relatively divisible} monomorphism (RD-monomorphism, for short) if it remains a monomorphism when tensored with cyclic modules of the form $R/Ra$, where $a \in R$ \cite[290]{wis}. A module is \emph{pure injective} (resp., \emph{coneat injective}, \emph{RD-injective}) if it is injective with respect to all pure (resp., coneat, RD-) monomorphisms \cite[p. 278]{wis} (\cite{cinj}). The classes $\mathcal{PI}$  of pure injective $R$-modules, $\mathcal{CI}$ of coneat injective $R$-modules, and, if $R$ is a commutative domain, $\mathcal{RDI}$ of RD-injective $R$-modules are all enveloping (\cite[p. 141]{eandj}, \cite{cinj}, and \cite[p. 425]{fs}). Moreover, $\Psi (\mathcal{PI})$, (resp., $\Psi (\mathcal{CI})$, $\Psi (\mathcal{RDI})$) consists of all pure (resp., coneat, RD-) monomorphisms and we have the injective structures $(\Psi (\mathcal{I}nj),\mathcal{I}nj)$, $(\Psi (\mathcal{PI}),\mathcal{PI})$, $(\Psi (\mathcal{CI})),\mathcal{CI})$, and $(\Psi (\mathcal{RDI}),\mathcal{RDI})$.
\begin{sloppypar}
%Recall that $\mathcal{M}$-\emph{copure} monomorphisms \cite[p. 322]{wis} are just the $\mathcal{F}$-monomorphisms. In this case, a submodule $A$ of a module $B$ is $\mathcal{M}$-\emph{copure} in $B$ if the inclusion $A \rightarrow B$ is an $\mathcal{F}$-monomorphism.
An $\mathcal{M}$-\emph{copure Baer injective} $R$-module \cite{cpbi} is one which is injective with respect to all $\mathcal{M}$-cpure monomorphisms $I\rightarrow R$. In \cite{rcinj}, it is proved that the class $\mathcal{C}BI$ of $\mathcal{M}$-copure Baer injective modules is preenveloping. So, $\Psi(\mathcal{C}BI)$ includes the $\mathcal{M}$-copure monomorphisms. It is easy to see that $\mathcal{C}BI = \Phi(\Psi(\mathcal{C}BI))$. Thus $(\Psi(\mathcal{C}BI), \mathcal{C}BI)$ is an injective structure.
\end{sloppypar}
One last example is when we have a hereditary torsion theory $\tau$ \cite{bland}. That is, a pair $(\mathcal{T},\mathcal{T}')$ of classes of modules such that both $\mathcal{T}$ and $\mathcal{T}'$ are maximal with the property that $\Hom(T,T')=0$ for every $T \in \mathcal{T}$ and $T' \in \mathcal{T}'$ and that $\mathcal{T}$ is closed for submodules. A submodule $A$ of a module $B$ is $\tau$-\emph{dense} in $B$ if $B/A\in \mathcal{T}$. A $\tau$-\emph{injective} module is one which is injective with respect to all $\tau$-dense monomorphisms. The class $\mathcal{\tau I}nj$ of $\tau$-injective modules is enveloping. ($\Psi(\mathcal{\tau I}nj)$, $\mathcal{\tau I}nj$) is an injective structure.
%\begin{definition} An injective structure $\mathcal{G}=(\mathcal{F},\mathcal{M})$ is called \emph{direct limit closed} if the class $\mathcal{F}$ is closed for direct limits.
%\end{definition}
In this paper, we are interested in injective structures whose classes of homomorphisms are closed for direct limits. The classes $\Psi (\mathcal{I}nj)$. $\Psi (\mathcal{PI})$, $\Psi (\mathcal{CI})$, and $\Psi (\mathcal{RDI})$ are all closed for direct limits, because direct limits commute with tensor products. Also, $\Psi (\mathcal{\tau I}nj)$ is closed for direct limits because $\Psi (\mathcal{\tau I}nj)$ consists of $\tau$-dense monomorphisms. On the other hand, $\Psi$($R$-Mod), consisting of all split monomorphisms, is not closed for direct limits. 

The class of homomorphisms being closed for direct limits gives the class of injective modules an extra property:
\begin{lemma} \label{env}
    If $(\mathcal{F},\mathcal{M})$ is  an injective structure such that $\mathcal{F}$ is closed for direct limits, then the class $\mathcal{M}$ is enveloping.
    \begin{proof}
        Let $A$ be a module and let $A \rightarrow M_i$, $i\in I$ be any well-ordered inductive system of $\mathcal{M}$-preenvelopes of $A$. Now, all the maps $A \rightarrow M_i$ must be in $\mathcal{F}$ because every module in $\mathcal{M}$ is injective with respect to them. But the injective structure $(\mathcal{F},\mathcal{M})$ is limit-closed, therefore, $A \rightarrow \varinjlim M_i$ is also in $\mathcal{F}$. So, any $\mathcal{M}$-preenvelope $A \rightarrow M$ of $A$ is injective with respect to the map $A \rightarrow \varinjlim M_i$. Hence, there is a factorization $A \rightarrow \varinjlim M_i \rightarrow M$ of $A \rightarrow M$. By \cite[lemma 6.6.1]{eandj}, $\mathcal{M}$ is enveloping.
    \end{proof}
\end{lemma}

\begin{lemma} \label{split} \cite[proposition 1]{maranda}
A module $M$ is $\mathcal{F}$-injective if and only if any $\mathcal{F}$-map $M \rightarrow N$ splits. In particular, such maps must be monomrphisms. 
\end{lemma}
%\begin{corollary} \label{dsummand} If a module $M$ is $\mathcal{F}$-injective then any monomorphism $M \rightarrow N$ in $\mathcal{F}$ splits.
%\end{corollary}
\section{\bf{Isomorphisms of Modules in $\mathcal{M}$}}
%Let $(\mathcal{F},\mathcal{M})$ be an injective structure on $R$-Mod such that the class $\mathcal{F}$ is closed for direct limits.
\begin{theorem} \label{main} %Suppose that there are two $\mathcal{F}$-maps $A \rightarrow B$ and $B \rightarrow A$ between two $\mathcal{F}$-injective modules $A$ and $B$.
If $(\mathcal{F},\mathcal{M})$ is an injective structure of $R$-Mod such that $\mathcal{F}$ is closed for direct limits, then any two $\mathcal{F}$-injective modules $A$ and $B$ with $\mathcal{F}$-homomorphisms $A \rightarrow B$ and $B \rightarrow A$
are isomorphic.
\begin{proof} 
As all $\mathcal{F}$-maps between $\mathcal{F}$-injective modules are split monomorphisms (lemma \ref{split}), we may assume that $B$ is a direct summand of $A$ and $f: A \rightarrow B$ a (split) $\mathcal{F}$-monomorphism. %Then, by lemma \ref{split}, both $f$ and $B \hookrightarrow A$ split.
So $A=B \oplus B'$ for some submodule $B'$ of $A$ and 
$$A=B \oplus B' \supseteq B' \oplus f(A) = B' \oplus f(B') \oplus f(B) \supseteq B' \oplus f(B')\oplus f^2(B')=\cdots$$
Therefore, $A\supseteq B' \oplus f(B')\oplus f^2(B') \oplus \cdots = C$. Now, each of the inclusions $f(B') \hookrightarrow B$, $f^2(B') \hookrightarrow B$, $\cdots$ is an $\mathcal{F}$-monomorphism. By assumption, their direct limit $f(B') \oplus f^2(B') \oplus \cdots = f(C) \hookrightarrow B$ is also an $\mathcal{F}$-monomorphism with respect to which any $\mathcal{F}$-injective module is injective. This means that $f(C) \hookrightarrow B$ is an $\mathcal{M}$-preenvelope of $f(C)$. %Moreover, $f(C)=f(B')\oplus f^2(B') \oplus \cdots$ is a direct summamd of $C$. Hence, $f(C) \hookrightarrow A$ is also an $\mathcal{M}$-preenvelope of $f(C)$.
If $M$ is an $\mathcal{M}$-envelope of $f(C)$ (Lemma \ref{env}) then by \cite[proposition 6.1.2]{eandj}, $M$ is a direct summand of $B$. Put $B=M \oplus M'$ so that $A= B' \oplus M \oplus M'$. Therefore, $B' \oplus M$ is $\mathcal{M}$-injective and by \cite[proposition 6.4.1]{eandj}, it is the injective envelope of $C = B'\oplus f(C)$. Now, as $f$ is a monomorphism, $C$ and  $f(C)$ must be isomorphic and therefore, so are their respective $\mathcal{M}$-envelopes $M$ and $B' \oplus M$. Thus, $B=M \oplus M' \cong B' \oplus M \oplus M'=A$.
\end{proof}
\end{theorem}
\begin{corollary}
    \begin{enumerate}
        \item Two (pure, coneat) injective modules isomorphic to (pure, coneat) submodules of each other are isomorphic. (cf. \cite[Theorem]{bumby}, \cite[Theorem 4]{diaz1}, and \cite[Theorem 11]{diaz2}.)
        \item For a hereditary torsion theory $\tau$, two $\tau$-injective modules isomorphic to $\tau$-dense submodules of each other are isomorphic.
        \item For a domain $R$, two RD-injective modules isomorphic to RD-submodules of each other are isomorphic. (cf. \cite[Theorem 4]{diaz1}.)
    \end{enumerate}
\end{corollary}
Theorem \ref{main} can be extended to include any two modules, not just $\mathcal{F}$-injective ones.
\begin{proposition}
   Let $(\mathcal{F},\mathcal{M})$ be an injective structure of $R$-Mod such that $\mathcal{F}$ is closed for direct limits. If $A$ and $B$ are any modules with $\mathcal{F}$-maps  $A \rightarrow B$ and $B \rightarrow A$ between them, then the $\mathcal{M}$-envelopes of $A$ and $B$ are isomorphic.
   \begin{proof}
       Let $A\rightarrow \mathcal{A}$ and $B \rightarrow \mathcal{B}$ be the $\mathcal{M}$-envleopes of $A$ and $B$, respectively. Let $f: A \rightarrow B$ and $g: B \rightarrow A$ be given $\mathcal{F}$-homomorphisms. Now, $g: B \rightarrow A$ and the $\mathcal{M}$-envelope $A\rightarrow \mathcal{A}$ are in $\mathcal{F}$ and, hence so is their composition $B\rightarrow \mathcal{A}$ with respect to which any module in $\mathcal{M}$ is injective. This means that $B \rightarrow \mathcal{A}$ is an $\mathcal{M}$-preenvelope of $B$. Therefore, by \cite[proposition 6.1.2]{eandj}, $\mathcal{B}$ is a direct summand of $\mathcal{A}$. Similarly, $\mathcal{A}$ is a direct summand of $\mathcal{B}$. So by Theorem \ref{main}, $\mathcal{A} \cong \mathcal{B}$.
   \end{proof}
\end{proposition}

\begin{corollary}
    \begin{enumerate}
        \item If $A$ and $B$ are two modules with two (pure, coneat) monomorphisms $A \rightarrow B$ and $B \rightarrow A$, then the (pure, coneat) injective envelopes of $A$ and $B$ are isomorphic. (cf. \cite[corollary 1]{bumby}, \cite[corollary 5]{diaz1}, and \cite[corollary 12]{diaz2}.)
        \item Over a domain $R$, if $A$ and $B$ are two modules with two  RD-monomorphisms $A \rightarrow B$ and $B \rightarrow A$, then the RD-injective envelopes of $A$ and $B$ are isomorphic. (cf. \cite[corollary 5]{diaz1}.)
    \end{enumerate}
\end{corollary}
%\begin{definition}
 %   Let $(\mathcal{F},\mathcal{M})$ be an injective structure such that the class $\mathcal{M}$ is enveloping. A module $A$ is called $\mathcal{M}$-\emph{quasi injective} if for any endomorphism $f$ of $\mathcal{M}(A)$, $f\phi (A) \subseteq A$.
%\end{definition}
%\begin{examples}
 %    \begin{enumerate}
            %\item $\mathcal{I}nj$-quasi injective modules are the usual quasi injectives.
            %\item $\mathcal{PI}$-quasi injective modules are the pure quasi injectives.
            %\item $\mathcal{CI}$-quasi injective modules are the coneat quasi injectives.
%        \end{enumerate}
 %   \end{examples}

%Suppose now that the $\mathcal{M}$-envelopes are monomphic. Then a module is $\mathcal{M}$-quasi injective if and only if it is a fully invariant submodule of its $\mathcal{M}$-envelope.

\bibliographystyle{amsplain}
%\bibliographystyle{apacite}
%\bibliography{mybib.bib}

\end{document}